\documentstyle[12pt]{article}
\begin{document}
\newtheorem{lem}{Lemma}
\newtheorem{rem}{Remark}
\newtheorem{question}{Question}
\newtheorem{prop}{Proposition}
\newtheorem{cor}{Corollary}
\newtheorem{thm}{Theorem}
\title
{A geometric problem and the Hopf Lemma. II
}
\author{YanYan Li\thanks{Partially
 supported by
      NSF grant DMS-0401118.}
\\
Department of Mathematics\\
Rutgers University\\
110 Frelinghuysen Road\\
Piscataway, NJ 08854\\
\\
Louis Nirenberg
 \\
Courant Institute\\
251 Mercer Street \\
New York, NY 10012\\
}

\date{ }
\maketitle

%%%%%  The following two lines should be changed
%\input amstex
\input { amssym.def}
%%%%%%%%%


\bigskip

\centerline{Dedicated to the memory of S.S. Chern}

\bigskip

\begin{abstract}
A classical result of A.D. Alexandrov states that a
connected compact smooth $n-$dimensional manifold without 
boundary, embedded in $\Bbb R^{n+1}$, and such that its mean curvature 
is constant, is a sphere.  Here we study the problem of symmetry of $M$ in a hyperplane $X_{n+1}=$constant
in case $M$ satisfies: for any
two points $(X', X_{n+1})$, $(X', \widehat X_{n+1})$ on
$M$, with $X_{n+1}>\widehat X_{n+1}$, the mean
curvature at the first is not greater than that at the 
second.   Symmetry need not always
hold, but in this paper, we
establish it under some additional conditions.  
Some variations of the Hopf Lemma are also presented.
Several open problems are described.
Part I dealt with corresponding
one dimensional problems.

\end{abstract}

\setcounter {section} {0}

\section{Introduction}

\subsection{}
In this  sequel to \cite{LN}, we continue our
study on a geometric problem related
to a classical result of A.D. Alexandrov.  
Throughout the paper $M$ is  a smooth compact connected 
embedded hypersurface in $\Bbb R^{n+1}$, 
its mean curvature is
$$
H(X):=\frac 1n\left[k_1(X)+\cdots +k_n(X)\right],
$$
where
$k(X)=(k_1(X), \cdots, k_n(X))$ denote the principle curvatures
of $M$ at $X$ with respect to the inner normal.
Let $G$ denote the open bounded set bounded by $M$. 

The problem we consider is to prove a symmetry property 
for $M$ satisfying the following

\medskip

\noindent{\bf Main Assumption.}\
For any two points $ (X', X_{n+1}), (X', \widehat X_{n+1})\in M$
satisfying $X_{n+1}\ge \widehat X_{n+1}$ and that
$\{(X', \theta X_{n+1}+(1-\theta)\widehat X_{n+1})\
|\ 0\le \theta\le 1\}$ lies in $\overline G$, 
\begin{equation}
H(X', X_{n+1})\le H(X', \widehat X_{n+1})
 \label{d1}
 \end{equation}
holds.

\medskip

It is suggested that the reader first read the introduction of
\cite{LN}.

If the mean curvature function $H$ is constant on $M$, then 
$M$ must be a standard sphere by a classical result in
 \cite{A}.   Under the Main Assumption and assuming that the
mean curvature function  can be extended to $\Bbb R^{n+1}$ as
a monotone Lipschitz function,  it was
proved in \cite{Li} that $M$ must be symmetric about
some hyperplane $X_{n+1}=constant$.
Examples given in \cite{LN} show that 
the Main Assumption alone is not enough to guarantee
the symmetry.  It is not difficult to see that
 the examples can be made
so  that  the mean curvature
function can be extended to a monotone function in
$\Bbb R^{n+1}$ which is  H\"older
continuous with H\"older exponent  $\alpha$ for
any $0<\alpha<1$.     The examples do not
satisfy

\medskip

\noindent{\bf Condition S.}\ $M$ stays on one side of any hyperplane
parallel to the $X_{n+1}$-axis that is tangent to $M$.

\medskip

\begin{rem} It is not difficult to see that
Condition S implies that $G$, the interior of
$M$, is convex in the $X_{n+1}$ direction.
The converse is not true.
\end{rem}

We make the following

\medskip

\noindent {\bf Conjecture.}\  
Any smooth compact connected
 embedded hypersurface $M$ in $\Bbb R^{n+1}$
 satisfying
the Main Assumption and Condition S must be symmetric about
some hyperplane $X_{n+1}=constant$.

\medskip

The Conjecture in dimension $n=1$ was proved in
\cite{LN}. 
A  crucial ingredient in the
proof was later established in \cite{B} by
a simpler method.
In the present paper we present 
results concerning the Conjecture in
dimensions $n\ge 2$.  
%We also present results 
%on related issues concerning
%variations of the Hopf  Lemma and the strong maximum
%principle.  These issues were naturally led to
%by our approach to the Conjecture and were
%studied in dimension $n=1$
% in \cite{LN}.
Our main result 
in the present paper for higher dimensions requires a further condition:

\medskip

\noindent{\bf Condition T.}\  Any line parallel to the 
$X_{n+1}$-axis that is tangent to $M$ has contact of finite order.

\medskip

If $\nu(X)=(\nu_1(X), \cdots, \nu_{n+1}(X))$ denotes the inner unit normal
of $M$ at $X$, we will consider the set 
$$
T:=\{X\in M\ |\ \nu_{n+1}(X)=0\}
$$
i.e. the set of points on $M$ where the tangent planes
are parallel  to
the $X_{n+1}$-axis.

For a point $\bar X$ in $T$, we often work in a new coordinate system which is
orthogonal to the original one.  The new coordinate system is 
centered at $\bar X$ and denoted by
$(y_1, \cdots, y_{n-1}, t, y_{n+1})$, with $y_{n+1}$-axis pointing 
in the direction of the inner normal of $M$ at $\bar X$, 
$t-$axis pointing to the opposite direction of the
$X_{n+1}$-axis, and the  $(y_1, \cdots, y_{n-1}, t)-$ coordinate
plane is the tangent plane of $M$ at $\bar X$.
In this new coordinate system,  let
 $v=v(t, y)$, $y=(y_1, \cdots, y_{n-1})$, denote the 
smooth function whose graph is $M$ near the origin.
Clearly $v(0,0)=0$ and $\nabla v(0,0)=0$.
With this notation, Condition T means:
For any $\bar X\in T$, there exists some integer $k\ge 2$ such that
\begin{equation}
\partial ^{k}_tv(0,0)\ne 0.
\label{g8}
\end{equation}

Our main theorem, Theorem \ref{thm2} below, also assumes that 
 $M$ is locally convex in the $X_{n+1}$-direction near $T$ in the following 
sense: 

\medskip

\noindent{\bf Condition LC.}\
 For every point $\bar X$ in $T$, if we view   $M$ locally
as the graph of a function defined on the tangent plane,
the function is convex in the $X_{n+1}$ direction near the
point.   Namely, the above defined function $v$ satisfies $v_{tt}\ge 0$
near the origin for every $\bar X$ in $T$.

\medskip

\begin{rem} Neither of the Condition S and LC implies the
other.
\end{rem}

Here is our main result.
\begin{thm} Let $M$  
 satisfy the Main Assumption and Conditions T and LC.
Then
 $M$ must be symmetric with respect to some hyperplane
 $X_{n+1}=constant$.
\label{thm2}
\end{thm}

\begin{cor}  Let $M$ be a  smooth 
compact convex hypersurface in $\Bbb R^{n+1}$
satisfying the Main Assumption and Condition T.
Then
 $M$ must be symmetric with respect to some hyperplane
 $X_{n+1}=constant$.  
\label{cor1}
\end{cor}

In particular, we have
\begin{cor}  Let $M$ be a real analytic compact convex hypersurface
in $\Bbb R^{n+1}$
satisfying the Main Assumption.
Then
 $M$ must be symmetric with respect to some hyperplane
 $X_{n+1}=constant$.  
\label{cor2}
\end{cor}

 The conclusion of Theorem \ref{thm2} still
holds when  the mean curvature function is replaced
by more general curvature functions.
Let $M$ 
 satisfy  Conditions T and LC, and let    $g(k_1, k_2, \cdots, k_n)$ be a $C^3$
function, symmetric in $(k_1, \cdots, k_n)$, defined
in an open neighborhood $\Gamma$ of
$
\{(k_1(X), \cdots, \cdots(X))\ |\
X\in M\}
$,  and satisfying in $\Gamma$
$$
\frac{\partial g}{\partial k_i}>0,\qquad 
\ \ 1\le i\le n
$$
and
$$
\frac{  \partial^2g}
{\partial k_i\partial k_j}\eta^i\eta^j\le 0,\qquad\forall\ \eta\in \Bbb R^n.
$$

\begin{thm} Let $M$ and  $g$ be as above.  We assume that
 for any two points $ (X', X_{n+1})$, 
$(X', \widehat X_{n+1})\in M$
  satisfying $X_{n+1}\ge \widehat X_{n+1}$
with  $\{(X', \theta X_{n+1} +(1-\theta)\widehat X_{n+1})\
|\ 0\le \theta\le 1\}$ lying inside $M$, 
 \begin{equation}
g(k(X', X_{n+1}))\le g(k(X', \widehat X_{n+1}))
\label{d1new}
 \end{equation}
holds.
  Then  $M$ must be symmetric with respect to some hyperplane
   $X_{n+1}=constant$.
\label{thm3}
\end{thm}

In \cite{LN} we mentioned extension by A. Ros of Alexandrov$'$s
result to other symmetric functions of the principal curvatures.
There was earlier work \cite{H} by P. Hartman.

Elementary symmetric functions satisfy the above properties of $g$ in
appropriate regions:
For $1\le m\le n$, let
$$
\sigma_m(k_1, \cdots, k_n)=
\sum_{ 1\le i_1<\cdots <i_m\le n}
k_{i_1}\cdots k_{i_m}
$$
be the $m-$th elementary 
symmetric functions,
and let
$$
g_m:= (\sigma_m)^{\frac 1m}.
$$
It is well known that
$g_m$ satisfies the above properties in 
$$
\Gamma_m:=\{(k_1, \cdots, k_n)\in \Bbb R^n\ |\
\sigma_j(k_1, \cdots, k_n)>0\
\mbox{for}\ 1\le j\le m\}.
$$

\medskip

\subsection{}
Theorem \ref{thm2} is proved in Section 2; our method of proof begins as in 
that of A.D. Alexandrov, using the method of moving planes.  As indicated in
\cite{LN} one is led to the need for extensions of the
classical Hopf Lemma.  Here we also present
some variations of the Hopf Lemma and the strong maximum 
principle.  In \cite{LN} these were studied in one dimension.

The Hopf Lemma is a local result.  We have not been able to prove the 
analogous local result for our problem.  Our proof of Theorem 
\ref{thm2}, which uses the maximum principle,
is via a global argument.

Here are some plausible variations of the
Hopf Lemma adapted for our problem.
Consider 
\begin{equation}
\Omega=\{(t,y)\ |\ y\in \Bbb R^{n-1}, |y|<1,
0<t<1\},
\label{D1-1new}
\end{equation}
\begin{equation}
u, v\in C^\infty(\overline \Omega),
\label{D1-2}
\end{equation}
\begin{equation}
u\ge v\ge 0,\qquad \mbox{in}\ \Omega,
\label{D6-1}
\end{equation}
\begin{equation}
u(0,y)=v(0,y),\quad \forall\ |y|<1;
\qquad u(0,0)=v(0,0)=0,
\label{D6-2}
\end{equation}
\begin{equation}
u_t(0,0)=0,
\label{D6-5}
\end{equation}
\begin{equation}
u_t>0,\qquad \mbox{in}\ \Omega,
\label{D6-3}
\end{equation}
and
\begin{equation}
\left\{
\begin{array}{l}
\mbox{whenever}\ u(t,y)=v(s,y), 0<s<1, |y|<1,\ \mbox{then there}\\
H(\nabla u, \nabla^2 u)(t,y)\le
H(\nabla v, \nabla^2 v)(s,y),
\end{array}
\right.
\label{D6-4}
\end{equation}
where
$$
H(\nabla u, \nabla^2 u)
:=\frac 1n\nabla \left( \frac{\nabla u}{ \sqrt{1+|\nabla u|^2} }\right)
$$
gives the mean curvature of the graph of $u$.

The followings are some plausible variations of the 
Hopf Lemma.

\medskip

\noindent{\bf Open Problem 1.}\
Assume the above.  Is it true that either
\begin{equation}
u\equiv v\ \ \mbox{near}\ (0,0)
\label{open1}
\end{equation}
or
\begin{equation}
v\equiv 0\ \ \mbox{near}\ (0,0)?
\label{open2}
\end{equation}

\medskip

A weaker version is

\medskip

\noindent{\bf Open Problem 2.}\ In addition to the assumption
in Open Problem 1, we further assume that
\begin{equation}
w(t,y):=
\left\{
\begin{array}{ll}
v(t,y),& t\ge 0, |y|<1\\
u(-t, y),&  t<0, |y|<1
\end{array}
\right.
\ \mbox{is}\ C^\infty\ \mbox{in}\
\{(t,y)\ |\ |t|<1, |y|<1\}.
\label{jjj}
\end{equation}
Is it true that either (\ref{open1}) or (\ref{open2}) holds?

\medskip

If the answer to Open Problem 2 is affirmative, then the Conjecture can
be proved 
by modification of the arguments in \cite{LN} and the
present paper.  The answer to Open Problem 1 
in dimension $n=1$ is affirmative, as proved
in \cite{LN}.  On the other hand, the answer to Open Problem 2
in higher dimensions is not known even under an additional hypothesis that
$\displaystyle{
\frac{\partial ^k v}{\partial t^k}(0,0)>0
}$ for some integer $k\ge 2$.  

Though our knowledge about the problems above concerning
variations of the Hopf Lemma is very limited, here is
a simple variation
 of the strong maximum
principle.

\begin{thm} For $n\ge 2$, let $\Omega$ be in (\ref{D1-1new}), and
let $u, v\in C^2(\Omega)$ satisfy (\ref{D6-4}),
\begin{equation}
u\ge v\qquad \mbox{in}\ \Omega,
\label{bb1}
\end{equation}
and
\begin{equation}
\max\{u_t, v_t\}>0\qquad \mbox{in}\ \Omega.
\label{bb2}
\end{equation}
Then either
\begin{equation}
u>v\qquad \mbox{in}\ \Omega,
\label{bb3}
\end{equation}
or
\begin{equation}
u\equiv v\qquad \mbox{in}\ \Omega.
\label{bb4}
\end{equation}
\label{thm5}
\end{thm}

A more general result, Theorem \ref{thm6}, is
proved in Section 4.

\begin{rem}  The analogue of Theorem \ref{thm5} in dimension $n=1$ was proved
in \cite{LN}. The same conclusion of Theorem \ref{thm5}
 holds when  the mean
curvature operator $H(\nabla u, \nabla^2 u)$ is
replaced by any
elliptic operator $F(u, \nabla u, \nabla^2 u)$,
see Theorem \ref{thm6} in Section 4.
\end{rem}

Another weaker form of Open Problem 1 is

\noindent{\bf Open Problem 3.}\ 
Let $u$ and $v$ satisfy (\ref{D1-2}), 
(\ref{D6-5}), (\ref{D6-3}), (\ref{D6-4}),
\begin{equation}
u\ge v>0\qquad\mbox{in}\ \Omega,
\label{ccc}
\end{equation}
and
\begin{equation}
u(0, y)=v(0, y)=0, \qquad \forall\ |y|<1.
\label{cccc}
\end{equation}
Is it true that (\ref{open1}) holds?

\medskip

In Open Problem 3, one may also replace the mean curvature operator by
other elliptic operators including the Laplacian operator.
In Section 5 we give some partial results concerning 
some of these open problems.

Theorems  \ref{thm2} and \ref{thm3} are proved in 
Section 2 and 3.  Section 5 contains some partial results on the
open problems 1-3 and variations
of the Hopf Lemma.
We think that they are of independent interest.

\section{Proof of Theorem \ref{thm2}}

This is the main section of the paper.

\subsection{} We start with the method of moving planes.

\noindent{ \bf Proof of Theorem \ref{thm2}.}\
Without loss of generality, we assume that
$$
\max\{X_{n+1}\ |\ (X_1, \cdots, X_{n+1})\in M
\ \mbox{for some}\ X_1, \cdots, X_n\}=0.
$$
For $\lambda<0$, let 
$S_\lambda:=
\{X\in M\ |\ X_{n+1}>\lambda\}$ denote the
portion of $M$ above the hyperplane $\{X_{n+1}=\lambda\}$,
and $S_\lambda'$ denote the mirror image of $S_\lambda$ with respect to 
$\{X_{n+1}=\lambda\}$. 
It is obvious that for $\lambda<0$ but close to $0$,
\begin{equation}
S_\lambda'\ \mbox{lies in}\
G, \mbox{the interior of}\ M, \ 
S_\lambda'\cap M=\emptyset,
\label{d2}
\end{equation}
and for all $X\in \partial S_\lambda'$,
\begin{equation}
\nu_{n+1}(X)<0.
\label{d3}
\end{equation}

Let $(\lambda_0, 0)$ denote the largest open interval
such that (\ref{d2}) and (\ref{d3}) hold for all $\lambda\in
(\lambda_0, 0)$.
To prove the theorem we need only  to show that
\begin{equation}
M=S_{\lambda_0}\cup \overline {  S_{\lambda_0}'},
\label{d4}
\end{equation}
where $\overline {  S_{\lambda_0}'}=  S_{\lambda_0}'
\cup \partial  S_{\lambda_0}'$.

It is easy to see from the definition of $\lambda_0$ that
\begin{equation}
\nu_{n+1}(X)<0,\qquad \forall\ X\in S_{\lambda_0},
\label{aaa}
\end{equation}
and that at least one of the following two cases occurs:
\begin{equation}
\mbox{there
exists some}\ \tilde X\in
S_{\lambda_0}'\cap M \ \mbox{with}\ \nu_{n+1}(\tilde X)>0,
\label{d5}
\end{equation}
\begin{equation}
\mbox{there
exists some}\ \tilde X\in 
\partial S_{\lambda_0}'\cap M\ \mbox{with}\
\nu_{n+1}(\tilde X)=0.
\label{d6}
\end{equation}

If (\ref{d5}) occurs,
$S_{\lambda_0}'$ and $M$ near $\tilde X$ can be represented as graphs of
smooth functions $u$ and $v$:
$$
u=u(X_1, \cdots, X_n),\ v=v(X_1, \cdots, X_n),
\quad \mbox{for}\ (X_1, \cdots, X_n)\ \mbox{close to}\
(\tilde X_1, \cdots, \tilde X_n).
$$
Clearly
\begin{equation}
u(\tilde X_1, \cdots, \tilde X_n)=v(\tilde X_1, \cdots, \tilde X_n),
\ \mbox{and}\ u\ge v\ \mbox{near}\ (\tilde X_1, \cdots, \tilde X_n).
\label{d7}
\end{equation}

By the Main Assumption,
$$
\nabla \left( \frac{\nabla u}{ \sqrt{1+|\nabla u|^2} }\right)
\le \nabla \left( \frac{\nabla v}{ \sqrt{1+|\nabla v|^2} }\right)
\quad\mbox{near}\
(\tilde X_1, \cdots, \tilde X_n).
$$
It follows, using the mean value theorem,  that 
$$
L(u-v):=a_{ij}\partial _{ij}(u-v)+b_i\partial_i(u-v)\le 0
\ \mbox{near}\ (\tilde X_1, \cdots, \tilde X_n),
$$
where $(a_{ij})$ is some smooth positive definite $n\times n$ matrix function
 and $\{b_i\}$ are some smooth 
functions, both near $(\tilde X_1, \cdots, \tilde X_n)$. 
By the strong maximum principle, in view of
(\ref{d7}), $u\equiv v$ near $(\tilde X_1, \cdots, \tilde X_n)$.
Since the argument applies to every point $\tilde X$ satisfying
(\ref{d5}), we obtain 
(\ref{d4}) in this case.

\medskip

Now we treat the much more delicate case (\ref{d6}), and we will
assume below that (\ref{d5}) does not occur.
Consider
\begin{equation}
M_{\lambda_0}:=\{X=(X_1, \cdots, X_{n+1})\in M\
|\ X_{n+1}<\lambda_0\},
\label{g30}
\end{equation}
 the part of $M$ below the hyperplane $T_{\lambda_0}$.
For $X\in M_{\lambda_0}$, let
\begin{equation}
\overline Y(X):= (X_1, \cdots, X_n, \lambda_0),
\label{d9}
\end{equation}
and
\begin{equation}
O:=\{X\in M_{\lambda_0}\
|\ \{\mbox{the segment between}\
X \mbox{and}\ \overline Y(X)\}\cap 
S_{\lambda_0}'\ne \emptyset\}.
\label{d10}
\end{equation}
For $X\in O$, we define
\begin{equation}
Y(X):=\{\mbox{
the segment between}\ X\ \mbox{and}\ \overline Y(X)\}\cap 
S_{\lambda_0'}.
\label{d11}
\end{equation}
It is clear that $Y(X)$ is uniquely defined
for $X\in O$, and it is a smooth function on 
$O$.

\medskip

***************************

{\bf Insert Figure 1
}

***************************

\medskip

Let
\begin{equation}
\tau(X):= \mbox{dist}(X, Y(X)),\quad
\bar \tau(X):= \mbox{dist}(X, \overline Y(X))
=\lambda_0- X_{n+1},\qquad X\in O,
\label{d12}
\end{equation}
where dist denotes the Euclidean distance between the two points.
Both $\tau(X)$ and $\bar\tau(X)$ are smooth functions on
$O$, and they can be extended continuously
to the closure of $O$.

Since we have assumed that (\ref{d5}) does not occur, 
\begin{equation}
\tau(X)>0\qquad \forall\ X\in O.
\label{f1}
\end{equation}

  Clearly,
\begin{equation}
\tau(X)=\bar \tau(X)\qquad\forall\ X\in \partial O.
\label{1}
\end{equation}

\subsection{}
The main step in our proof of Theorem \ref{thm2} is to establish
\begin{prop} Assume (\ref{f1}).  
Then  there exist some constants $\epsilon, c>0$ such that
\begin{equation}
\tau(X)\ge c \bar \tau(X), \qquad \forall\ X\in O_\epsilon:=\{X\in O\ 
|\ 
\bar \tau(X)=\lambda_0-X_{n+1}< \epsilon\}.
\label{d13}
\end{equation}
\label{prop1}
\end{prop}
\begin{rem} Proposition \ref{prop1} holds without 
assuming Condition T.
\end{rem}

\noindent {\bf Proof of Proposition \ref{prop1}.}\
If $O_\epsilon=\emptyset$ for
some $\epsilon>0$, (\ref{d13}) is considered to hold trivially.
So we assume that $O_\epsilon\ne \emptyset$ for all $\epsilon>0$.
In fact, Condition T guarantees that  $O_\epsilon$ is not empty,
as shown towards the end of the proof of Theorem \ref{thm2}.
For small $\epsilon>0$, by (\ref{f1}) and (\ref{1}),
there exists some small number $c=c(\epsilon)\in (0, \frac 18)$, 
depending on $\epsilon$,
such that
\begin{equation}
\tau(X)\ge 4c\bar\tau(X)\ge 2c\left(
\bar \tau(X)+\bar\tau(X)^{\frac 32}\right),\quad
\forall\ X\in \partial O_\epsilon.
\label{f2}
\end{equation}
For sufficiently small $\epsilon$ we will prove (\ref{d13}) with
$c=c(\epsilon)$ arguing by contradiction.

Here is a sketch of how the argument goes.

From our Main Assumption (\ref{d1}) it follows that $\tau$ satisfies a second 
order linear differential inequality on $O_\epsilon$, though we
do not write it down at a general point.  If (\ref{f2}) fails,
there is a point $\tilde X$ where
$$
\sigma:=
\tau-2c\left(\bar\tau +\bar \tau^{\frac 32}\right)
$$
has a negative minimum in $O_\epsilon$.  But  
\begin{equation}
\lim_{\epsilon\to 0} \sup_{x\in O_\epsilon}
dist\left(x, T\cap \{X_{n+1}=\lambda_0\}\right)=0.
\label{O1}
\end{equation}
Thus $\tilde X$ has a closest point 
 $\bar X$ in
 $\left(T\cap\{X_{n+1}=\lambda_0\}\right)$.
We then use our special coordinates, taking
$\bar X$ as origin, and compute, near $\bar X$, the differential
inequality.  $L\tau<0$.
In addition, we find that
$$
L\left(\bar \tau+\bar \tau^{\frac 32}\right)>0.
$$
But then, $\tilde X$ cannot be a minimum point
for $\sigma$.

We now proceed with the argument.  First we write down
the linear inequality
$L\tau<0$ near any point
$\bar X$ in
 $\left(T\cap\{X_{n+1}=\lambda_0\}\right)$, working
in the new coordinate system
\newline 
$(y_1,\cdots, y_{n-1}, t, y_{n+1})$ as 
described earlier.

Let, with $y=(y_1, \cdots, y_{n-1})$ and $\delta>0$ 
some small universal number,
$$
\Omega=\{(t,y)\ |\ 0<t<\delta,  y\in \Bbb R^{n-1}, |y|<\delta\},
$$
and
\begin{equation}
\Omega^+=\{(s,y)\in \Omega\ |\ 
\mbox{there exists some } 0<t<s\ \mbox{such that}\
u(t,y)=v(s,y)\}.
\label{ggg}
\end{equation}
Throughout the paper, a number is said to be universal if
it depends only on $M$.
We note that $(s,y)\in \Omega^+$
if and only if $(s,y,v(s,y))$ lies in $O$.

By (\ref{f1}) and (\ref{aaa}),
\begin{equation}
u(t, y):= v(-t, y)\qquad \mbox{in}\ \Omega
\label{r3}
\end{equation}
satisfies
$$
u(t, y)>v(t,y),\qquad (t,y)\in \Omega,
$$
and
\begin{equation}
u_t(t,y)>0,\qquad (t,y)\in \Omega.
\label{g1}
\end{equation}
With (\ref{g1}), an application of the implicit function theorem
yields that for any $(s,y)\in \Omega^+$, there exists a unique
$t=t(s,y)\in (0, s)$ satisfying
\begin{equation}
u(t(s,y))=v(s,y),
\label{g2}
\end{equation}
and the function $t(s,y)$ is smooth in $\Omega^+$.

By the Main Assumption,
\begin{equation}
H(\nabla u, \nabla^2 u)(t(s,y), y)\le 
H(\nabla v, \nabla^2 v)(s,y)\qquad \forall\ (s,y)\in \Omega^+.
\label{g3}
\end{equation}

Set
\begin{equation}
\tau(s,y)=s-t(s,y), \qquad (s,y)\in \Omega^+.
\label{g40}
\end{equation}
Differentiating (\ref{g2}), we have, with $1\le \alpha, \beta\le n-1$,
\begin{eqnarray}
v_s(s,y)&=&u_t(t,y)-u_t(t,y)\tau_s(s,y),
\label{341}\\
v_{y_\alpha}(s,y)&=& u_{y_\alpha}(t,y)-u_t(t,y)\tau_{y_\alpha}(s,y),\label{342}\\
v_{ss}(s,y)&=& u_{tt}(t,y)-\tau_s(s,y)[2-\tau_s(s,y)]
u_{tt}(t,y)-u_t(t,y)\tau_{ss}(s,y),\label{343}\\
v_{sy_\alpha}(s,y)&=& v_{y_\alpha s}(s,y)=
u_{ty_\alpha}(t,y)-u_{tt}(t,y)\tau_{y_\alpha}(s,y)[1-\tau(s,y)]
\nonumber\\
&& \qquad \qquad
\qquad-u_{ty_\alpha}(t,y)\tau_s(s,y)-u_t(t,y)\tau_{sy_\alpha}(s,y),\label{344}\\
v_{y_\alpha y_\beta}(s,y)&=&
u_{y_\alpha y_\beta}(t,y)
-u_{ty_\alpha}(t,y)\tau_{y_\beta}(s,y) -u_{ty_\beta}(t,y)\tau_{y_\alpha}(s,y)
\nonumber\\
&&
+u_{tt}(t,y)\tau_{y_\alpha}(s,y)\tau_{y_\beta}(s,y)-
u_t(t,y)\tau_{y_\alpha y_\beta}(s,y).\label{345}
\end{eqnarray}

By the mean value theorem, we have, with
$t=t(s,y)$ and $(s,y)\in \Omega^+$ and using
(\ref{341}) and (\ref{342}),
\begin{eqnarray}
&&H(\nabla v(s,y), \nabla^2 v(s,y))-H(\nabla u(t,y), \nabla^2 v(s,y))
\nonumber\\
&=& \left(\int_0^1 H_p(\theta \nabla v(s,y)+(1-\theta) \nabla u(t,y),
 \nabla^2 v(s,y)) d\theta\right)
\cdot \left( \nabla v(s,y) -\nabla u(t,y) \right)\nonumber\\
&=& \left[O(1) \tau_s(s,y)+
O(1)\cdot \nabla_y\tau (s,y)\right]  u_t(t,y),
\label{4-1}
\end{eqnarray}
where $O(1)$ satisfies $|O(1)|\le C$ for some universal constant
$C$.

Next
we have, using (\ref{343}), (\ref{344})
and (\ref{345}), 
\begin{eqnarray}
&&H(\nabla u(t,y), \nabla^2 v(s,y))-
H(\nabla u(t,y),  \nabla^2 u(t,y))\nonumber\\
&=& -u_t(t,y)\left[H_{00}\tau_{ss}(s,y)
+2\sum_{\alpha=1}^{n-1}
H_{0\alpha}\tau_{sy_\alpha}
+ \sum_{\alpha, \beta=1}^{n-1} H_{\alpha\beta} \tau_{y_\alpha y_\beta}
\right]\nonumber\\
&&-H_{00}(2-\tau_s)u_{tt}\tau_s -2
\sum_{\alpha=1}^{n-1} H_{0\alpha}u_{t y_\alpha}\tau _s
+\eta\cdot \nabla _y \tau,
\label{5-1}
\end{eqnarray}
where $H_{ij}$ denotes
$\displaystyle{
\frac{  \partial H(\nabla u(t,y), N)  }
{\partial N_{ij}   }
}$ which are
 independent of the matrix $N$, and $\eta$ denotes some vector-valued
function in $L^\infty_{loc}(\Omega^+)$ which may vary from line
to line.
Note that here and in the following,
$\nabla u$ denotes $\nabla u(t,y)$, etc.,
$\nabla v$ denotes $\nabla v(s,y)$, $H_{00}$
denotes $\displaystyle{
\frac{\partial H}{\partial u_{tt}}
}$  etc.

We deduce from (\ref{g3}), (\ref{4-1}) and (\ref{5-1})
that
\begin{eqnarray}
0&\le& 
-u_t\left[ H_{00}\tau_{ss}
+2 \sum_{\alpha=1}^{n-1} H_{0\alpha} \tau_{sy_\alpha}
+ \sum_{\alpha, \beta=1}^{n-1} H_{\alpha\beta} \tau_{y_\alpha y_\beta}
\right]\nonumber
\\
&& -H_{00} (2-\tau_s) u_{tt}\tau_s
+O(1) u_t \tau_s -2  \sum_{\alpha=1}^{n-1} H_{0\alpha} u_{ty_\alpha}
\tau_s +\eta\cdot \nabla_y \tau.
\label{6-1}
\end{eqnarray}

Since
$$
H_{0\alpha}=-(1+|\nabla u|^2)^{ -\frac 32} u_t u_{y_\alpha},
\qquad 1\le \alpha\le n-1,
$$
we have
\begin{eqnarray}
0&\le&
-u_t\left[ H_{00}\tau_{ss}
+ 2\sum_{\alpha=1}^{n-1} H_{0\alpha} \tau_{sy_\alpha}
+ \sum_{\alpha, \beta=1}^{n-1} H_{\alpha\beta} \tau_{y_\alpha y_\beta}
\right]\nonumber
\\
&& -H_{00} (2-\tau_s) u_{tt}\tau_s
+O(1) u_t \tau_s+\eta\cdot \nabla_y \tau.
\label{7-1}
\end{eqnarray}

Define
\begin{equation}
L:=  H_{00}\partial_{ss}
+2 \sum_{\alpha=1}^{n-1} H_{0\alpha}\partial_{sy_\alpha}
+ \sum_{\alpha, \beta=1}^{n-1} H_{\alpha\beta} \partial_{y_\alpha y_\beta}
+ H_{00} (2-\tau_s) \frac{   u_{tt} }{ u_t }
\partial_s -O(1)\partial_s -\frac {\eta_\alpha }
{ u_t} \partial_{y_\alpha}.
\label{w1}
\end{equation}
We know from (\ref{7-1}) that
\begin{equation}
L\tau\le 0\qquad\mbox{in}\ \Omega^+.
\label{w2}
\end{equation}

Let
\begin{equation}
\hat \tau(s,y):= s+s^{1+\bar \epsilon},
\label{8-1}
\end{equation}
where $\bar \epsilon=\frac 12$.

A calculation gives
$$
L\hat \tau=
H_{00} (2-\tau_s) \frac { u_{tt} }{ u_t}
\frac{d}{ds} [s+s^{ 1+\bar \epsilon }]
+H_{00}  (1+\bar \epsilon) \bar \epsilon s^{ \bar \epsilon-1}
-O(1)[1+(1+\bar \epsilon) s^{ \bar \epsilon}].
$$

        \begin{lem} There exists some universal constant
$\delta'>0$ such that
$$
\tau_s(s,y)<1,\qquad \forall\
(s,y)\in\Omega^+, \ |(s,y)|<\delta'.
$$
	        \label{lem1new}
		        \end{lem}

			        \noindent{\bf Proof.}\
In view of (\ref{341}) and the positivity of $u_t$ in
$\Omega^+$, we only need to show that $v_s(s,y)>0$.
We prove this by contradiction.
Suppose that  $v_s(s, y)=0$ for some
small $(s,y)$, $s>0$.
Recall that $u(t,y)=v(s,y)$, $0<t=t(s,y)<s$.
Since $M$ satisfies Condition LC,
\begin{equation}
v(0,y)\ge v(s,y)+v_s(s,y)(-s)=v(s,y).
\label{nnn}
\end{equation}
It follows that
$$
u(0, y)=v(0, y)\ge v(s,y)=u(t, y),
$$
which  violates $u_t>0$ in $\Omega^+$.

\vskip 5pt
\hfill $\Box$
\vskip 5pt

We will assume from now on, making $\delta$ smaller if necessary,
that $\delta\le \delta'$.
Since $M$ satisfies Condition LC,
 we have,
making $\delta$ smaller if necessary,
\begin{equation}
u_{tt}\ge 0,\qquad \mbox{in}\ \Omega^+.
\label{w3}
\end{equation}
It follows, using (\ref{g1}), (\ref{w3}) and
Lemma \ref{lem1new}, that
$$
L\hat \tau\ge  (1+\bar \epsilon) \bar \epsilon H_{ N_{00} } 
s^{ \bar \epsilon-1} +O(1).
$$
Thus,  making 
$\delta$ smaller if necessary,
\begin{equation}
L\hat \tau >0 \ \mbox{in}\ \Omega^+.
\label{w4}
\end{equation}

Now the value of $\delta$ is fixed; it works
works for every $\bar X$ in
$T\cap \{X_{n+1}=\lambda_0\}$.  We see from 
(\ref{O1}) that for $\epsilon>0$ small,
\begin{equation}\sup_{x\in O_\epsilon}
dist\left(x, T\cap\{X_{n+1}=\lambda_0\}\right)<\frac \delta 2.
\label{A2-1}
\end{equation}
As we described above, we fix such an $\epsilon$ now and
take $c=c(\epsilon)$ the one
in (\ref{f2}).  We will prove (\ref{d13}) 
with this value of $c$ arguing by contradiction.  Suppose
that (\ref{d13}) does not hold, then there exists
$\tilde X\in \overline O_\epsilon$ such that
\begin{equation}
\left(\tau-c(\bar\tau+\bar\tau^{\frac 32})\right)(\tilde X)
=\min_{ \overline O_\epsilon }
\left(\tau-c(\bar\tau+\bar\tau^{\frac 32})\right)<0.
\label{B2-1}
\end{equation}
Because of (\ref{f2}),
$\tilde X\in O_\epsilon$.
Namely, $\tilde X$ is an interior local minimum point
of $\tau-c(\bar\tau+\bar\tau^{\frac 32})$ in $O_\epsilon$.
Let $\bar X$ be a closest point
in $T\cap \{X_{n+1}=\lambda_0\}$ to
$\tilde X$.  We know from (\ref{A2-1}) that 
$dist(\bar X, \tilde X)\le \frac \delta 2$.  With this $\bar X$
and the function $v$ and $u$ defined earlier, $\tilde X$ 
corresponds to some $(\tilde s, \tilde y)$
in $\Omega^+$, with $0<\tilde s<\epsilon$.
Clearly $(\tilde s, \tilde y)$ is an interior local
maximum point of $\tau-c\hat\tau$ in
$\Omega^+$.  Thus
$$
L(\tau-c\hat\tau)\ge 0\ \ 
\mbox{at}\ (\tilde s, \tilde y).
$$
On the other hand, by (\ref{w2}) and (\ref{w4}),
$$
L(\tau-c\hat \tau)< 0\ \mbox{at}\ (\bar s, \bar y).
$$
A contradiction.  
Proposition \ref{prop1}
is established.

\vskip 5pt
\hfill $\Box$
\vskip 5pt

Now we use Condition T to show that
(\ref{d6}) cannot hold if (\ref{d5}) does not occur.
 
Since we are treating case (\ref{d6}), 
we let 
$\bar X$ be a point satisfying (\ref{d6}),
$v=v(t,y)$ be the function defined earlier 
for the point, and, in view of Condition T,
let $k\ge 2$ be the smallest $k$ 
satisfying (\ref{g8}). 
Set $u(t,y):=v(-t,y)$ for $t\ge 0$.
By the definition of $\lambda_0$
and by the assumption that case (\ref{d5}) does not occur,
$u(t,y)>v(t,y)$ for $(t,y)$ small and $t>0$.
Since $M$ satisfies Condition LC,
 $v_{tt}(t,0)\ge 0$ 
for small $t$.
 So
 $k$ is even, 
$ \partial^k_t v(0,0)>0$ and therefore $v(t, 0)>0$ for small
$t>0$ which clearly implies that
$(t, 0)\in \Omega^+$ for small $t>0$.

Now
$$
v(t,0)= \frac 1{k!}  \partial^k_t v(0,0) t^k+O(t^{k+1}),
\quad u(t,y)= \frac 1{k!}  \partial^k_t v(0,0) t^k+O(t^{k+1}).
$$
From
$
u(t(s,0), 0)=v(s,0)
$ and the above, we see easily that
\begin{equation}
\displaystyle{
\lim_{s\to 0^+}\frac {t(s,0) }s=1}.
\label{g11}
\end{equation}
Since
$$
\bar \tau(s, 0, v(s,0))
=s, \qquad
\tau(s, 0, v(s,0)) =\tau(s,0)=s-t(s,0),
$$
and $(s, 0, v(s,0))\in O$, we know from Proposition \ref{prop1} that
for some constant $c>0$ and for all $s>0$ small,
$$
s-t(s,0)=\tau(s, 0, v(s,0)) \ge c \bar \tau(s, 0, v(s,0))=cs.
$$
This is contradicted by (\ref{g11}).  Theorem \ref{thm2} is established.

\vskip 5pt
\hfill $\Box$
\vskip 5pt

\section{Proof of Theorem \ref{thm3}}

\noindent{ \bf Proof of Theorem \ref{thm3}.}\
We follow the proof of Theorem \ref{thm2} until (\ref{d7}).
Let $A(\nabla u, \nabla^2 u):=(A_{il}(\nabla u, \nabla^2 u)$ 
denote the second fundamental form of 
$S_{\lambda_0}'$ with respect to its first fundamental form. 
Then, see lemma 1.1 of \cite{CNS}, 
$$
A_{il}(\nabla u, \nabla ^2 u)=
\frac 1 w
\left\{ u_{il}-
\frac{  u_i u_j u_{jl} }
{ w(1+w) }
-\frac {  u_l u_k u_{ki} }
{  w(1+w) }
+ \frac {   u_i u_l u_j u_k u_{jk}  }
{  w^2 (1+w)^2  } \right\},
$$
where
$$
w=\sqrt{  1+|\nabla u|^2 }.
$$

Let
 $ {\cal S}^{n\times n}$ denote the set of
real symmetric $n\times n$ matrices, and let
$O(n)$ denote the set of $n\times n$ 
real orthogonal matrices.  For $A\in {\cal S}^{n\times n}$ we use
$k(A)$ to denote $(k_1(A), \cdots, k_n(A))$ where
$ k_1(A), \cdots, k_n(A)$ are the $n$ eigenvalues of $A$.
We define a function $G$ on
$$
U:=\{A\in {\cal S}^{n\times n}\ |\ k(A)\in \Gamma\}
$$
by 
$$
G(A):=g(k(A)).
$$
By the properties of $g$,  $G\in C^3(U)$,
 $$
 O^{-1}UO=U\quad\forall\ O\in O(n),
 $$
\begin{equation}
G_{A_{ij}}(A)\eta^i\eta^j>0, \qquad \forall A\in U,\ \eta\in \Bbb R^n
\setminus\{0\},
\label{a3new}
\end{equation}
\begin{equation}
G(O^{-1}AO)=G(A)\quad\forall\ 
A\in U\ \mbox{and}\ O\in O(n),
\label{g15}
\end{equation}
\begin{equation}
G_{A_{ij}A_{kl}}(A)\xi^{ij}\xi^{kl}\le 0,\qquad
\forall\ A\in U, \forall\ \xi\in {\cal S}^{n\times n}.
\label{a4new}
\end{equation}

By (\ref{d1new}),
$$
A(\nabla u, \nabla^2 u), A(\nabla v, \nabla ^2v)\in U
\ \mbox{near}\ (\tilde X_1, \cdots, \tilde X_n),
$$
and
$$
G(A(\nabla u, \nabla^2 u))\le G(A(\nabla v, \nabla ^2v)).
$$
Using the mean  value theorem as
usual we have, by
 (\ref{a3new}), 
$$
L(u-v):=a_{ij}\partial _{ij}(u-v)+b_i\partial_i(u-v)\le 0
\ \mbox{near}\ (\tilde X_1, \cdots, \tilde X_n),
$$
where $(a_{ij})$ is some smooth positive definite $n\times n$ 
matrix function
 and $\{b_i\}$ are some smooth
 functions, both near $(\tilde X_1, \cdots, \tilde X_n)$.
We obtain
 (\ref{d4}) in this case as in the proof of Theorem \ref{thm2}.

Now we treat the much more delicate case (\ref{d6}), and we will
assume below that (\ref{d5}) does not occur.
We follow from (\ref{g30}) until (\ref{1}), and we give the

\medskip

\noindent{\bf Proof of 
Proposition \ref{prop1} under the hypotheses of 
Theorem \ref{thm3}.}\
Follow from the beginning of the proof of Proposition \ref{prop1}
until (\ref{g2}). Instead of (\ref{g3}), we have
\begin{equation}
F(\nabla u,
\nabla ^2u)(t(s,y), y)\le F(\nabla v, \nabla^2 v)(s,y),
\label{g3new}
\end{equation}
where we have used notation
$$
F(\nabla u, \nabla ^2 u)=
G(A(\nabla u, \nabla ^2 u)).
$$

With $\tau(s,y)$ defined in (\ref{g40}), we still have 
(\ref{341})-(\ref{345}).
Similar to (\ref{4-1}), we have 
\begin{eqnarray}
&&F(\nabla v(s,y), \nabla^2 v(s,y))-F(\nabla u(t,y), \nabla^2 v(s,y))
\nonumber\\
&=& \left[O(1) \tau_s(s,y)+
O(1)\cdot \nabla_y\tau (s,y)\right]  u_t(t,y).
\label{D1-1}
\end{eqnarray}

Since we only work in regions where $u_t$ and $\epsilon$ are very small,
there
$$
\{A(\nabla u(t,y), \theta \nabla^2 v(s,y)
+(1-\theta)\nabla^2 u(t,y))\ |
\ 0\le \theta\le 1\}\subset U.
$$

Since $G(A)$ is concave in $A$ and $A(p,N)$ is linear in $N$,
$F(p,N)$ is concave in $N$.  So we have
\begin{eqnarray*}
&&
F(\nabla u(t,y), \nabla^2 v(s,y))-
F(\nabla u(t,y), \nabla^2 u(t,y))\\
&\le &
F_{jk}(\nabla u(t,y), \nabla^2 u(t,y))\bigg[v_{jk}(s,y)-u_{jk}(t,y)\bigg],
\end{eqnarray*}
where 
$$
F_{jk}(p,N):= \frac{ \partial F(p,N) }{
\partial N_{jk}  }.
$$

It is easy to see from
(\ref{a3new})
  that  for some universal constant 
$C_1>1$, 
$$
\frac 1{C_1}|\xi|^2\le F_{jk}(\nabla u(t,y), \nabla^2 u(t,y))\xi_j\xi_k
\le C_1 |\xi|^2,\qquad \forall\
\xi\in \Bbb R^n.
$$

Next, still with $t=t(s,y)$ and  $(s,y)\in \Omega^+$
and using (\ref{343}), (\ref{344}) and (\ref{345}),
we have
\begin{eqnarray}
&&F(\nabla u(t,y), \nabla^2 v(s,y))-
F(\nabla u(t,y),  \nabla^2 u(t,y))\nonumber\\
&\le& -u_t(t,y)\left[F_{00}\tau_{ss}(s,y)
+2\sum_{\alpha=1}^{n-1}
F_{0\alpha}\tau_{sy_\alpha}
+ \sum_{\alpha, \beta=1}^{n-1} F_{\alpha\beta }\tau_{y_\alpha y_\beta}
\right]\nonumber\\
&&-F_{00}(2-\tau_s)u_{tt}\tau_s -2
\sum_{\alpha=1}^{n-1} F_{0\alpha}u_{0 \alpha}\tau _s
+\eta\cdot \nabla _y \tau,
\label{D3-1}
\end{eqnarray}
where $u_{0\alpha}=u_{t y_\alpha}$, $F_{jk}$ denotes
$F_{jk}(\nabla u(t,y), \nabla^2 u(t,y))$, and 
 $\eta$ denotes some vector-valued
function in $L^\infty_{loc}(\Omega^+)$ which may vary from line
to line.

The term $-F_{00}(2-\tau_s)u_{tt}\tau_s$
can be handled as in the proof of Theorem \ref{thm2}, by
using  Lemma \ref{lem1new} and Condition LC.
  We mainly need to
show that
\begin{equation}
\sum_{\alpha=1}^{n-1} F_{0\alpha}u_{0\alpha}=O(1)u_t.
\label{main1}
\end{equation}

For $1\le \alpha\le n-1$,
$$
F_{0\alpha}u_{0\alpha}=
\sum_{0\le i,l\le n-1}G_{il}\cdot \frac {\partial  A(\nabla u,N)_{il} }
		     {\partial N_{0\alpha} } u_{0\alpha}.
		       $$

  Observe that
  \begin{equation}
   \frac {\partial A(\nabla u, N)_{il} }
    {\partial N_{0\beta} }=
     \frac 1w\left(\delta_{i0}\delta_{l\beta}+O(1) u_t\right),\qquad
      1\le  \beta\le n-1,
         \label{a14newnew}
	   \end{equation}
	     It follows that
	       \begin{eqnarray}
	         \sum_{\alpha=1}^{n-1}  F_{0\alpha}\cdot u_{0\alpha}
		   &=&  \sum_{\alpha\ge 1}\sum_{l\ge 0}  G_{0l}
		     \frac {\partial  A(\nabla u,N)_{0l} }
		     {\partial N_{0\alpha} } u_{0\alpha}+O(1)u_t\nonumber\\
   &=&  \sum_{\alpha,  \beta\ge 1}  G_{0\beta}
		     \frac {\partial  A(\nabla u,N)_{0\beta} }
		     {\partial N_{0\alpha} } u_{0\alpha}+O(1)u_t.
		       \label{a19newnew}
		         \end{eqnarray}
Since $\{A(\nabla u, N)\}_{\beta\ge 1}$ is linear in
$\{N_{0\alpha}\}_{\alpha\ge 1}$, we have, using  Lemma
\ref{lema-1}, 
\begin{eqnarray*}
  \sum_{\alpha\ge 1}  F_{0\alpha}\cdot u_{0\alpha}
		   &=&   \sum_{\beta\ge 1}G_{0\beta}(A)A_{0\beta}
+O(1)u_t=O(1)\sum_{\beta\ge 1}|A_{0\beta}|^2+O(1)u_t\\
&=&O(1)\sum_{\beta\ge 1}|u_{0\beta}|^2+O(1)u_t.
\end{eqnarray*}

			     Since $u_t\ge 0$
			     and $\nabla^2_y u_t=O(1)$
			     in $2\Omega$,
			      we have,
			     using some well known inequality,
		      see  \cite{NT},  for some universal
			       constant $C$,
			      \begin{equation}
			      \sum_{j\ge 1}|u_{ty_j}(t,y)|\le C
			      \sqrt{u_t(t,y)}\qquad
			      \forall\ (t,y)\in  \Omega.
			      \label{a27new}
			      \end{equation}
		      With this, we obtain (\ref{main1}).

We deduce from (\ref{g3new}), (\ref{D1-1}),
(\ref{D3-1}) and (\ref{main1}) that
\begin{equation}
L\tau \le 0,\qquad\mbox{in}\
\Omega^+
\label{C1-1}
\end{equation}
where
$$
L:=F_{00}\partial_{ss}+2\sum_{\alpha=1}^{n-1}
F_{0\alpha}\partial_{s y_\alpha}
+\sum_{\alpha, \beta=1}^{n-1} F_{\alpha\beta}
\partial_{ y_\alpha y_\beta}
+F_{00} (2-\tau_s)\frac{t_{tt}}{u_t}\partial _s
-O(1)\partial_s
-\frac{\eta_\alpha}{u_t} \partial_{y_\alpha}.
$$

Using Condition S and Lemma \ref{lem1new}, as
in the proof of Theorem \ref{thm2}, we have
\begin{equation}
F_{00} (2-\tau_s) \frac{u_{tt}}{u_t}\ge 0,
\qquad \mbox{in}\ \Omega.
\label{C2-1}
\end{equation}

Let $\hat \tau(s,y)$ be in (\ref{8-1}) with
$\bar\epsilon=\frac 12$, we derive from (\ref{C2-1}) that
for $\delta>0$ small,
\begin{equation}
L\hat\tau>0,\qquad\mbox{in}\ \Omega^+.
\label{C2-2}
\end{equation}

With (\ref{C1-1}) and (\ref{C2-2}), the rest of the proof
of Theorem \ref{thm3} follows as in the proof of Theorem \ref{thm2}.

\vskip 5pt
\hfill $\Box$
\vskip 5pt

\section{A variation of the strong maximum principle}

In this section we establish
a result more general than
Theorem \ref{thm5}.
We consider $F\in C^1(\Bbb R\times \Bbb R^n\times {\cal S}^{n\times n})$
satisfying
$$
\frac{\partial F}{\partial N_{ij}}(s,p,N)\xi_i\xi_j>0,
\qquad\forall\ \xi\in \Bbb R^n\setminus\{0\}, \ \forall\ (s,p,N).
$$

\begin{thm} For $n\ge 2$, let $F$ be as above, and
 let $\Omega$ be in (\ref{D1-1new}).  We assume that 
 $u, v\in C^2(\Omega)$ satisfy (\ref{bb1}), (\ref{bb2})
and
\begin{equation}
\left\{
\begin{array}{l}
\mbox{if}\ u(t,y)=v(s,y), 0<s<1, |y|<1,\ \mbox{then there}\\
F(u, \nabla u, \nabla^2 u)(t,y)\le
F(v, \nabla v, \nabla^2 v)(s,y).
\end{array}
\right.
\label{D6-4new}
\end{equation}
Then either (\ref{bb3}) or (\ref{bb4}) holds.
\label{thm6}
\end{thm}

\begin{rem}  The analogue of Theorem \ref{thm6} in dimension $n=1$ was proved
in \cite{LN}.
\end{rem}

\noindent{\bf Proof.}\ Suppose that (\ref{bb3}) 
does not hold, then $u(\bar s, \bar y)=v(\bar s, \bar y)$
for some $(\bar s, \bar y)\in \Omega$.  Clearly, $u_t(\bar s, \bar y)=
v_t(\bar s, \bar y)>0$ and, by the implicit function theorem,
for $(s,y)$ close to $(\bar s, \bar y)$ there exists a unique $C^2$ function
$t=t(s,y)$ such that
$u(t(s,y),y)=v(s,y)$.  Thus
$F(u,\nabla u, \nabla^2 u)(t(s,y), y)
\le F(v, \nabla v, \nabla^2v)(s,y)$.
As in the proof of Theorem \ref{thm2},
(\ref{341})-(\ref{345}) hold near 
$(\bar s, \bar y)$ with $\tau(s,y)=s-t(s,y)$.
As usual these lead to 
$
L\tau\le 0\qquad \mbox{near}\ (\bar s, \bar y)
$ 
where
$L=a_{ij}\partial_{ij}+b_i\partial$ with 
$(a_{ij})$ positive definite.  Since $\tau(\bar s, \bar y)=0$ and
$\tau\ge 0$ near $(\bar s, \bar y)$, we have, by the
strong maximum principle, $\tau\equiv 0$ near
$(\bar s, \bar y)$.  Namely $u\equiv v$ near
$(\bar s, \bar y)$.  Theorem
\ref{thm6} is established.

\vskip 5pt
\hfill $\Box$
\vskip 5pt

\section{Partial results on Open Problems 1-3}

In this section we give some partial results on or
related to Open Problems 1-3 and variations of the Hopf Lemma.

\subsection{}
\begin{thm}  Let $\Omega$ be as in (\ref{D1-1new}), and let
$u$ and $v$ satisfy (\ref{D1-2}), (\ref{D6-2}),  (\ref{D6-5}), (\ref{ccc}),
(\ref{D6-3}),
(\ref{D6-4}) and
\begin{equation}
u_{tt}\ge 0\qquad\mbox{in}\ \Omega.
\label{mmm}
\end{equation}
  We assume that for some
open set $0\in \omega\subset \overline \omega\subset
 \{y\in \Bbb R^n\
|\ |y|<1\}$,
\begin{equation}
\frac{\partial v}{ \partial t}(0, y)<0, \qquad\forall\ y\in \partial\omega.
\label{eee}
\end{equation}
Then
$$
\frac{  \partial^ k u}
{\partial t^k}(0, 0)= 0,\qquad\forall\ k\ge 2.
$$
\label{thm7}
\end{thm}
\begin{rem} It is clear from the proof that the conclusion 
of Theorem \ref{thm7} still holds when  the
mean curvature operator is replaced by the more general curvature
operators in Theorem \ref{thm6}.
\end{rem}

\noindent{\bf Proof.}\  We prove it by contradiction argument.
Suppose the contrary, then
for some integer $k\ge 2$,
\begin{equation}
\frac{  \partial^ ku}
{\partial t^k}(0, 0)\ne 0,
\ \frac{  \partial^ iu}
{\partial t^i}(0, 0)=0,\ 1\le i\le k-1.
\label{fff}
\end{equation}
By (\ref{eee}), (\ref{D6-3}) and 
 Theorem \ref{thm5},  $u>v$ in $\Omega$.
 Let $\Omega^+$ be as in
(\ref{ggg}). 
 Clearly for some $\epsilon_1>0$,
$(t, 0)\in \Omega^+$ for all $0<t<\epsilon_1$.
On the other hand,  in view of (\ref{eee}) and  (\ref{D6-3}),
there exists some $\epsilon_2>0$ such that
$\displaystyle{
\{(t, y)\ |\ 0<t<\epsilon_2\}\cap \Omega^+
=\emptyset}$, $\forall\ y\in \partial \omega$.  
For $0<\epsilon<\min\{\epsilon_1, \epsilon_2\}$, let
$$
\Omega^+_\epsilon:=\{ (s,y)\in \Omega^+\ |\
0<s<\epsilon\}.
$$
Then $\Omega^+_\epsilon$ is 
a nonempty open set satisfying
\begin{equation}
\Omega^+_\epsilon\subset
\{(s,y)\ |\ 0<s<\epsilon, y\in \omega\}.
\label{kkk}
\end{equation}
Let $t(s,y)$ and $\tau(s,y)$ be defined as in
(\ref{g2}) and (\ref{g40}),  then
$\tau =s$ on $\partial \Omega^+_\epsilon\cap \{0<s<\epsilon\}$,
and $\tau>0$ on $\partial \Omega^+_\epsilon\cap \{s=\epsilon\}$.
Thus, for some constant $c=c(\epsilon)\in (0, \frac 14)$,
\begin{equation}
\tau-c(s+s^{\frac 32})\ge 0\qquad\mbox{on}\ \partial \Omega^+_\epsilon.
\label{hhh}
\end{equation}
Let $L$ be defined in (\ref{w1}). By (\ref{mmm}), we still
have (\ref{nnn}), and therefore we still have
$\tau_s<1$ in $\Omega^+$.
  Making $\epsilon$ smaller
if necessary, we have, as established in the
proof of Theorem \ref{thm2},
$
L\tau\le 0 <L (s+s^{ \frac 32})$ in $ \Omega^+_\epsilon.
$
Thus
$$
L\left(\tau-c(s+s^{\frac 32})\right)<0,\qquad \mbox{in}\ \Omega^+_\epsilon.
$$
It follows that $\tau\ge
c(s+s^{\frac 32})\ge cs$ in $\Omega^+_\epsilon$.
With (\ref{fff}), we reach a contradiction by using 
the argument towards the end of the proof of
Theorem \ref{thm2}.

\vskip 5pt
\hfill $\Box$
\vskip 5pt

\subsection{}
Let $\Omega$ be as in (\ref{D1-1new}), and let
\begin{equation}
f\in C^\infty([-1, 1]^{n-1}\times (0,\infty)),
\label{1-0a}
\end{equation}
\begin{equation}
u\in C^\infty(\overline \Omega),\ u>0\ 
\mbox{in}\ \Omega, 
\label{1-0b} \end{equation}
\begin{equation}
u(0, y)=0\qquad \forall\ |y|<1,
\label{1-0c}
\end{equation}
and
\begin{equation}
\Delta u(t, y)=f(y, u(t,y)), \qquad \mbox{in}\ \Omega.
\label{1-1}
\end{equation}

Assume, for some integer $k\ge 1$,
\begin{equation}
u(t, y)= t^k a_k(y)+O(t^{k+1}),
\label{1-2}
\end{equation}
where
\begin{equation}
a_k(y)>0\quad\forall\ |y|\le 1.
\label{1-3}
\end{equation}

\begin{thm}  Let $\Omega$ and $ f$ be
as above, and let $u$ be a solution of (\ref{1-1}) satisfying
(\ref{1-0b}), (\ref{1-0c}),  (\ref{1-2}) and (\ref{1-3}).

\noindent (i)\ If $k=1$, then all  $\displaystyle{
\left\{
\frac {\partial^l}{ \partial t^l} u(0,y)
\right\}_{l\ge 2}
}$
are
determined by $f$ and $a_1(y)$.

\noindent (ii)\ If $k\ge 2$, then both $k$ and  $\displaystyle{
\left\{
\frac {\partial^l}{ \partial t^l} u(0,y)
\right\}_{l\ge k}
}$
are
determined by $f$.

\noindent (iii)\ If both $u$ 
and $v$ are solutions of (\ref{1-1}) satisfying
(\ref{1-0b}), (\ref{1-0c}),  (\ref{1-2}) and (\ref{1-3}),
so that by (i) and (ii), 
\begin{equation}
\frac {\partial^l}{ \partial t^l} u(0,y)
= \frac {\partial^l}{ \partial t^l} v(0,y),\qquad
\forall\ |y|<1, l\ge k,
\label{k2}
\end{equation}
and $u\ge v$ in $\Omega$,
then
$v\equiv u$ in $\Omega$.
\label{thm4}
\end{thm}

\begin{rem}  The $f$ in Theorem \ref{thm4} is not assumed to be smooth
up to $u=0$, otherwise the conclusion follows from 
classical results.
\end{rem}

First 
\begin{lem} Assume (\ref{1-0a})-(\ref{1-3}) with 
$k\ge 2$.
Then, for some constant $C>0$,
\begin{equation}
\sup_{ |y|\le 1, 0<s<1}
|f(y,s)-a_k(y)^{ \frac 2k}
k(k-1) s^{ \frac {k-2}k }|
s^{ \frac {1-k}k}<\infty,
\label{3-1}
\end{equation}
and
\begin{equation}
\lim_{t\to 0^+}
\frac{  f(y, u(t,y))  }
{\left[
\frac{u(t,y)}{ t^2}\right]
}
=k(k-1),\quad
\mbox{uniform in}\ |y|\le 1.
\label{3-2}
\end{equation}
Consequently, both $k$ and $a_k(y)$ are determined by $f$.
\label{lemma2}
\end{lem}

\noindent{\bf Proof.}\
Write
\begin{equation}
u(t,y)=t^ka_k(y)
+O(t^{k+1}).
\label{3-3}
\end{equation}
Then
$$
\Delta u(t,y)=k(k-1)t^{k-2}a_k(y)+O(t^{k-1}).
$$
Set 
$$
s=u(t,y)= t^k a_k(y)+O(t^{k+1}).
$$
We have
$$
t=
\left[ \frac s{ a_k(y)} \right]^{ \frac 1k}
\left[ 1+O(s^{\frac 1k})\right],
$$
$$
t^{k-2} =
\left[ \frac s{ a_k(y)} \right]^{ \frac {k-2}k}
+O( s^{\frac 1k}),
$$
$$
\Delta u(t,y)=k(k-1)
 a_k(y)^{\frac 2k}
s^{ \frac{k-2}k}
+O(s^{ \frac {k-1}k }).
$$
Estimate (\ref{3-1}) follows from 
this  and (\ref{1-1}).  It is easy to see from
(\ref{3-1}) that $k$ is determined by $f$.  In turn, again from
(\ref{3-1}), $a_k(y)$ is determined by $f$.

By (\ref{3-1}), we have, for some constant $C>0$, 
\begin{equation}
|f(y, u(t,y))-k(k-1)
a_k(y)^{\frac 2k}
u(t, y)^{ \frac{k-2}k} |
\le C u(t,y) ^{  \frac {k-1}k },
\qquad \forall\ |y|\le 1.
\label{4-1a}
\end{equation}
By (\ref{3-3}),
\begin{equation}
|u(t,y)-
t^k a_k(y)|\le Ct^{k+1},
\qquad
\forall\ |y|\le 1, 0<t<1.
\label{5-1a}
\end{equation}
By (\ref{4-1a}) and (\ref{5-1a}),
\begin{equation}
\lim_{ t\to 0^+}
\frac {  f(y, u(t,y))  }
{   u(t, y)^{  \frac{k-2}k }  }
=  k(k-1) a_k(y)^{ \frac 2k},
\qquad \forall\ |y|\le 1,
\label{5-2}
\end{equation}
and
\begin{equation}
\lim_{ t\to 0^+}
\frac {  u(t,y)^{\frac 2k}  }
{t^2} =a_k(y) ^{ \frac 2k},
\qquad \forall\ |y|\le 1.
\label{5-3}
\end{equation}
Estimate (\ref{3-2}) follows from (\ref{5-2}) and (\ref{5-3}).
Lemma \ref{lemma2} is established.

\vskip 5pt
\hfill $\Box$
\vskip 5pt

\noindent{\bf Proof of Part (i) and (ii) of Theorem \ref{thm4}.}\ 
Because of Lemma \ref{lemma2}, we only need to prove that
 $\displaystyle{
\left\{
\frac {\partial^l}{ \partial t^l} u(0,y)
\right\}_{l\ge k+1}
}$
are 
determined by $f$ and $a_k(y)$.  We will prove it by induction.

Write
\begin{equation}
u(t, y)=
 t^k a_k(y)
+ t^{k+1} a_{k+1}(y)
+\cdots +
t^{m-1} a_{m-1}(y)
+t^m a_m(y)
+O(t^{m+1}),
\label{6-1a}
\end{equation}
and we assume that $m\ge k+1$, and
$a_k(y), \cdots, a_{m-1}(y)$ are determined by
$f$.  We will prove that $a_m(y)$ is also determined by $f$
and $a_k(y)$.

Let
$$
s:=u(t,y)=
 t^k a_k(y)
+ t^{k+1} a_{k+1}(y)
+\cdots +
+t^m a_m(y)
+O(t^{m+1}).
$$
Then
$$
\lambda:=\left[ \frac s {a_k(y)} \right]^{ \frac 1k}
=t \left\{
1+ t \frac{a_{k+1}}{a_k}+\cdots +t^{m-k}\frac{a_m}{a_k}
+O( t^{m-k+1})\right\}^{\frac 1k}.
$$
It follows that
$$
\lambda=t\left\{1+\cdots+ t^{m-k-1} b_{m-k-1}
+t^{m-k} \frac{a_m}{ ka_k}
+O(t^{m-k+1})\right\},
$$
where $\{b_i(y)\}_{ i\le m-k-1}$ are determined by $f$
and $a_k(y)$.

Clearly, $\displaystyle{
\lim_{t\to 0}\frac \lambda t=1 }$.
We know that
$\displaystyle{
\frac{d\lambda}{ dt}|_{t=0},
\frac{d^2\lambda}{dt^2}|_{t=0}, \cdots,
\frac {d^{m-k}\lambda}{dt^{m-k}}|_{t=0}
}$ are
determined by $f$ and $a_k(y)$.  We now write $t$ in terms of $\lambda$.
First
$$
\frac{dt}{d\lambda}  \frac{d\lambda}{dt}=1.
$$
Applying $\frac d{d\lambda}$ to the above $m-k+1$
times, we have
$$
\frac{d^2 t}{ d\lambda^2} 
 \frac{d\lambda}{dt}+(\frac{dt}{d\lambda})^2
 \frac{d^2 \lambda}{ d t^2}=0,
 $$
 $$
 \frac{d^3 t}{d\lambda^3}  \frac{d\lambda}{dt}+
 \cdots+(\frac{dt}{d\lambda})^3
  \frac{d^3 \lambda}{ d t^3}=0,
  $$
  $$
  \cdots\cdots,
  $$
  $$
   \frac{d^{m-k+1} t}{d\lambda^{m-k+1} }   \frac{d\lambda}{dt}+
    \cdots+(\frac{dt}{d\lambda})^{m-k+1}  \frac{d^{m-k+1} \lambda}{ d t^
    {m-k+1}}=0.
    $$
Set $\lambda=0$ in the above.  All the ``$\cdots$'' contribute to
quantities 
determined by $f$ and $a_k(y)$.  Therefore 
$\frac{dt}{d\lambda}|_{\lambda=0}, 
\frac{d^2t}{d\lambda^2}|_{\lambda=0},
\cdots, \frac{d^{m-k}t }{ d\lambda^{m-k}}|_{\lambda=0}$ are
 determined by $f$, and
$\displaystyle{
\frac{ d^{m-k+1}t}{ d\lambda^{m-k+1}}|_{\lambda=0}+
\frac{ d^{m-k+1}\lambda}{ dt^{m-k+1}}|_{t=0}
}
$ is determined by $f$ and $a_k(y)$.  We also note that
$\displaystyle{
\frac{ d^{m-k+1}\lambda}{ dt^{m-k+1}}|_{t=0}=
(m-k+1)! \frac {a_m(y)}{ka_k(y)}
}$.  It follows that
$$
t=\lambda+\lambda^2 c_2(y)+\cdots+\lambda^{m-k}c_{m-k}(y)
-\lambda^{m-k+1}\frac {a_m(y)}{ka_k(y)}
+\lambda^{m-k+1}c_{m-k+1}(y)+O(\lambda^{m-k+2}),
$$
where $c_2(y), \cdots, c_{m-k+1}(y)$ are determined by $f$
and $a_k(y)$.

Applying $\Delta$ to (\ref{6-1a}) yields
$$
\Delta u(t,y)=\sum_{j=k}^{m-1} j(j-1) t^{j-2}\alpha_j(y)
+ m(m-1) t^{m-2} a_m(y)+O(t^{m-1}),
$$
where $\{\alpha_j(y)\}_{ k\le j\le m-1}$ are determined by $f$.
Since $f(y,u)=\Delta u$, we have
$$
f(y,s)= \sum_{j=k}^{m-1} j(j-1) t^{j-2}\alpha_j(y)
+ m(m-1) t^{m-2} a_m(y)+O(t^{m-1}).
$$
For $k\le j\le m-1$,
\begin{eqnarray*}
t^{j-2}&=&\lambda^{j-2}
\left\{ 1+\lambda c_2+\cdots+\lambda^{m-k-1}
c_{m-k}-\lambda^{m-k}\frac {a_m(y)}{ka_k(y)}
+O(\lambda^{m-k+1})\right\}^{j-2}\\
&=& 
\lambda^{j-2}
\left\{ 1+\lambda d_2+\cdots+\lambda^{m-k-1}
d_{m-k}
-\lambda^{m-k} \frac{j-2}{ka_k(y)} a_m(y)+
O(\lambda^{m-k+1})\right\}\\
&=& \lambda^{j-2}+\lambda^{j-1}d_2+
\cdots+\lambda^{ m+j-k-3}d_{m-k}
-  \frac{j-2}{ka_k(y)} \lambda^{ m+j-k-2}a_m(y)+
O(\lambda^{m+j-k-1}),
\end{eqnarray*}
where $d_2, \cdots, d_{m-k}$ are determined by $f$
and $a_k(y)$.

The coefficient of $a_m(y)$ in the above expansion of
$t^{k-2}$ is  of order $\lambda^{m-2}\sim 
t^{m-2}$,
while the coefficients of $a_m(y)$ in the
expansions of $t^{j-2}$ for
$k<j\le m-1$ are of higher order.
Thus
\begin{eqnarray*}
f(y,s)&=&\lambda^{k-2}e_{k-2}(y)+\cdots
+\lambda^{m-3} e_{m-3}(y)-
(k-1)(k-2)\lambda^{m-2} a_m(y)\\
&&
+m(m-1) \lambda^{m-2}a_m(y)+
\lambda^{m-2}e_{m-2}(y)+O(\lambda^{m-1}),
\end{eqnarray*}
where $\{e_j(y)\}_{ k-2\le j\le m-2}$ are determined by $f$
and $a_k(y)$.
Since $m\ge k+1$, we have
$m(m-1)>(k-1)(k-2)$. Therefore $a_m(y)$ is
also determined by $f$ and $a_k(y)$.  Part (i) and
(ii) of Theorem \ref{thm4} are
 established.

\vskip 5pt
\hfill $\Box$
\vskip 5pt

To prove Part (iii) of Theorem \ref{thm4}, 
we can make use of the following
\begin{thm} 
Let $w\in C^\infty(\overline {\Omega})$ satisfy
$$
w\ge 0\qquad \mbox{in}\ \Omega,
$$
$$
\partial^\alpha w(0,y)=0\quad \forall\ |y|\le 1,\ \forall\ \alpha,
$$
and, for some positive constant $C_0$,
$$
\Delta w\le C_0\frac w{t^2},\qquad\mbox{in}\ \Omega.
$$
Then
$$
w\equiv 0\quad\mbox{in}\ \Omega.
$$
\label{thmb1}
\end{thm}

Theorem \ref{thmb1} is an immediate corollary of the following kind of
Hopf Lemma.

\begin{thm} Consider a domain $\Omega$ in
$\Bbb R^n$ with $C^2$ boundary, and a positive function
$w$ in $\Omega$, $w\in C^\infty(\overline\Omega)$,
satisfying:  for some positive
constant $C_0$,
\begin{equation}
\Delta w(x) \le C_0 \frac{w(x)}{   dist(x,\partial\Omega)^2 }.
\label{97}
\end{equation}
Suppose $w=0$ at some boundary point $P$.
Then, along the inner normal to $\partial \Omega$ 
at $P$, close to $P$,
\begin{equation}
w(x)\ge a|x-P|^k
\label{98}
\end{equation}
where $a$ is a positive constant and $k>n$ satisfies
\begin{equation}
k(k-n)=C_0.
\label{99}
\end{equation}
\label{thm7prime}
\end{thm}

\noindent{\bf Proof.}\  Let $B_R$ be an open ball in $\Omega$ whose
boundary touches $\partial \Omega$ only at $P$. 
We may suppose that its center is the origin and that
$$
P=(-R, 0, \cdots, 0).
$$
Set $|x|=r$.  By (\ref{97}), $w$ satisfies
\begin{equation}
\Delta w(x)\le C_0 \frac{w(x)}{ (R-r)^2}\quad \mbox{in}\ B_R.
\label{100}
\end{equation}
We construct a comparison function
$$
h=(R-r)^k,
$$
with $k$ satisfying (\ref{99}).  In the region
$$
K:=\{x\in B_R\ |\ x_1<-\frac R2\}
$$
we have
$$
\frac{R-r}r<1.
$$
Then, in $K$,
$$
\Delta h=(R-r)^{k-2}
\left[ k(k-1) -(n-1)k \frac {R-r}r\right]
\ge (R-r)^{k-2}
k(k-n).
$$
Thus
$$
\Delta h\ge C_0\frac h{  (R-r)^2}.
$$
Since $w>0$ in $\Omega$, on the straight part of $\partial K$,
$$
w\ge ch
$$
for some constant $c>0$.  This same
inequality holds on the curved part of $\partial K$ since, there, $h=0$.
By the maximum principle it follows that
$$
w\ge ch\quad \mbox{in}\ K,
$$
and so (\ref{98}) follows.

\vskip 5pt
\hfill $\Box$
\vskip 5pt

By choosing $K$ much narrower one sees that
(\ref{98}) holds provided
$k(k-1)>C_0$; of
course $a$ depends on $k$.
An immediate consequence of this is
the following kind of Hopf Lemma, in which we may
take $k<2$.

\begin{cor} In a domain $\Omega$ in
$\Bbb R^n$ with $C^2$ boundary,
let $w\ge 0$, $w\in
C^2(\overline \Omega)$, satisfy 
(\ref{97}) near $\partial\Omega$, with $C_0<2$.
Suppose that at some boundary point $P$, $w$ and
its normal derivative vanish.  Then $w\equiv 0$.
\end{cor}

\begin{rem}  The proof of Theorem \ref{thm7prime}
applies also to a function
$w>0$ satisfying an elliptic
inequality
$$
Lw(x)\le C_0
\frac{  w(x)}
{  dist(x, \partial\Omega)^2}.
$$
Here
$ Lw=a_{ij}w_{ij}+b_i w_i+cw$
is uniformly elliptic with bounded coefficients.
The value of $k$ is, of course,
different.
\end{rem}

Returning to Theorem \ref{thm4}, we derive 
some further properties of $f$.
\begin{lem} Assume (\ref{1-0a})-(\ref{1-3}) with 
$k=1$.  Then
\begin{equation}
\sup_{ 0<s_1<s_2<1, |y|<1}
\frac{|f(y, s_1)-f(y,s_2)| }{ |s_1-s_2|}<\infty.
\label{1-4}
\end{equation}
\label{lemma1}
\end{lem}

\noindent{\bf Proof.}\ Write
$$
u(t,y)=ta_1(y)+t^2 a_2(y)+t^3 a_3(y)+O(t^4).
$$
Then
$$
u_t(t, y)=a_1(y)+O(t),
\quad
\Delta u_t(t,y)=6a_3(y)+\Delta_y a_1(y)+O(t).
$$
Applying $\partial_t$ to  (\ref{1-1})
yields
$$
f_u(y, u(t,y))=\frac{ \Delta u_t(t,y) }{u_t(t,y)}
=\frac{  6a_3(y)+\Delta_y a_1(y)+O(t) }
{  a_1(y)+O(t)  }.
$$
Let
$$
s=u(t,y)=ta_1(y)+O(t^2).
$$
Then 
$$
t=\frac s{  a_1(y) }+O(s^2).
$$
It follows that
$$
f_u(y,s)=
\frac{  6a_3(y)+\Delta_y a_1(y)+O(s) }
{  a_1(y) +O(s) }.
$$
This implies (\ref{1-4}).  Lemma \ref{lemma1} is established.

\vskip 5pt
\hfill $\Box$
\vskip 5pt

\begin{lem} If $k\ge 2$, 
there exists some positive constant $C$ such that
$$
f_u(y, s)\le Cs^{-\frac 2k},\qquad
\forall\ |y|\le 1, 0<s<1.
$$
\label{lemma5}
\end{lem}

\noindent{\bf Proof.}\
For $k\ge 3$,
$$
u_t(t,y)=k a_k(y)t^{k-1}+O(t^k),
\quad
\Delta u_t(t,y)=k(k-1)(k-2)a_k(y)t^{k-3}+ O(t^{k-2}).
$$
Applying $\partial_t$ to (\ref{1-1}) gives
\begin{equation}
\Delta u_t(t,y)=f_u(y, u)u_t(t,y).
\label{n1}
\end{equation}
Write
$$
s=u=a_k(y)t^k+O(t^{k+1}),
$$
we have 
$$
t=\left[ \frac s{a_k(y)} \right]^{\frac 1k}[ 1+O(s^{\frac 1k})].
$$
Thus
\begin{eqnarray*}
s^{\frac 2k}f_u(y,s)&=&
s^{\frac 2k} \frac { \Delta u_t(t,y) }{ u_t(t,y) }
= \frac{  k(k-1) (k-2) t^{k-3} +O(t^{k-2}) }
{  kt^{k-1} +O(t^k)  }
\to (k-1)(k-2) \ \mbox{as}\ s\to 0^+.
\end{eqnarray*}

For $k=2$, 
write
$$
u(t,y)=a_2(y)t^2+a_3(y)t^3+O(t^4).
$$
Applying $\partial_t$ to the above gives
$$
u_t(t,y)=2a_2(y)t+O(t^2),\quad
\Delta u_t(t,y)=6a_3(y)+O(t).
$$
We still have
(\ref{n1}).  Write
$$
s=u(t,y)=a_2(y)t^2 +O(t^3),
$$
we have
\begin{eqnarray*}
s^{\frac 12} f_u(y,s)&=&
s^{\frac 12} \frac  { \Delta u_t(t,y) }{ u_t(t,y) }
= s^{\frac 12} \frac  { 6a_3(y) +O(t) }
{2ta_2(y)+O(t^2) }\to
\frac{ 3a_3(y)} {\sqrt{a_2(y)}}\ \mbox{as}\ s\to 0^+.
\end{eqnarray*}
Lemma \ref{lemma5} is established.

\vskip 5pt
\hfill $\Box$
\vskip 5pt

Now the

\medskip

\noindent{\bf Proof of Part (iii) of Theorem \ref{thm4}.}\
For $k=1$, it follows from Lemma \ref{lemma1} that
$$
\Delta (u-v)=O(1)(u-v)\qquad \mbox{in}\ \Omega.
$$
Since
$$
u-v\ge 0\ \mbox{in}\ \Omega, \qquad
\frac{\partial }{\partial t}(u-v)(0, y)
=0\ \ \forall\ |y|<1,
$$
we have, by the Hopf Lemma and the strong maximum principle,
that
$u\equiv v$ in $\Omega$.  

Now we assume that  $k\ge 2$. Clearly, for any $\epsilon\in (0,1)$,
there exists some positive constant $C$ such that
\begin{equation}
\frac 1C t^k\le u, v\le Ct^k\qquad \mbox{in}\
(1-\epsilon)\Omega.
\label{v1}
\end{equation}
Using the equation satisfied
by $u$ and $v$, Lemma \ref{lemma5}
and (\ref{v1}),
we have
\begin{eqnarray*}
\Delta (u-v)&=& f(y, u)-f(y, v)
=\int_0^1 f_u(y, \theta u+(1-\theta)v)d\theta (u-v)\\
&= & O(1)(u-v)\int_0^1[  \theta u+(1-\theta)v ]^{ -\frac 2k }  
= O(1) \frac {u-v}{ t^2},\qquad \mbox{in}\ (1-\epsilon)\Omega.
\end{eqnarray*}
If $u\ge v$ in $\Omega$, then, by Theorem \ref{thmb1}, $u\equiv v$ in $(1-\epsilon)\Omega$.
Part (iii) of Theorem \ref{thm4}, where we have
$u\ge v$ in $\Omega$, 
 is established.

\vskip 5pt
\hfill $\Box$
\vskip 5pt

\subsection{}
The following two theorems are not used in this
paper.
\begin{thm}  For $n=2$, let $\Omega$ be in (\ref{D1-1new}), and let
$u$ and $v$ satisfy (\ref{D1-2}), (\ref{D6-2}),  (\ref{D6-5}), (\ref{ccc}),
(\ref{D6-3}), (\ref{jjj}), (\ref{mmm}) and
\begin{equation}
\left\{
\begin{array}{l}
\mbox{if}\ u(t,y)=v(s,y), 0<s<1, |y|<1,\ \mbox{then there}\\
\Delta u(t,y)=
\Delta v(s,y),
\end{array}
\right.
\label{D6-4newnew}
\end{equation}
We also assume that if $u_t(0, \bar y)=0$ for
some $|\bar y|<1$, then for some integer $\bar k\ge 2$, 
which may depend
on $\bar y$, 
\begin{equation}
\frac{  \partial^ {\bar k}u}
{\partial t^{\bar k}}(0, \bar y)\ne 0.
\label{fffnew}
\end{equation}
Then $u\equiv v$ in $\Omega$.
\label{thm8}
\end{thm}

\noindent{\bf Proof.}\  It is easy to see that we only need to consider the 
following two cases.

\medskip

\noindent{\it Case 1.} \  There exist $-1< \alpha<\beta<1$
 such that $u_t(0, y)=0$ for all $y\in (\alpha, \beta)$.

\noindent{\it Case 2.} \ There exist $-1<y^-<0<y^+<1$ such that
$u_t(0, y^\pm)>0$.

\medskip

In Case 1, we can find some point $\bar y\in (\alpha, \beta)$ and some
even integer $k\ge 2$ such that
$\partial^k_t u(0, \bar y)=\partial^k_t v(0, \bar y)>0$, and
$\partial ^i_tu(0, y)=\partial ^i_t v(0,y)=0$ for all $1\le i\le k-1$ and
all $y$ in some neighborhood of $\bar y$.
Without loss of generality, $\bar y=0$. By subtracting 
$u(0, y)$ from both $u$ and $v$, we may assume without loss
of generality that (\ref{cccc}) holds.   Now
(\ref{1-2}) holds with $k!
a_k(y)=\partial^k_t u(0,  y)=\partial^k_t v(0,  y)$.
Thus, for some $\delta>0$,
\begin{equation}
u\ge v>0\qquad\mbox{in}\ \delta \Omega.
\label{27-1}
\end{equation}
By (\ref{D6-3}), the map $(t, y)\to (u(t,y), y)$ is
a local diffeomorphism and,  by the
implicit function theorem, $t$ is 
locally a smooth function of $u$ and $y$ in $\Omega$.
Thus, in view of (\ref{D6-3}) and (\ref{cccc}), 
$$
\Delta u=f(y, u)\qquad \mbox{in}\ \delta\Omega
$$
where $f$ is some unknown smooth function in
$\{(y,u)\ |\ u>0, |y|<1\}$ and continuous 
in $\{(y,u)\ |\ u\ge 0, |y|<1\}$.
By (\ref{27-1}),
for every $(s,y)\in \delta\Omega$, there exists 
some $(t,y)$, with $0<t<s$, such that
$u(t,y)=v(s,y)$.  Thus, by (\ref{D6-4newnew}),
$$
\Delta v=f(y,v)\qquad\mbox{in}\ \delta\Omega.
$$
  An application of Theorem \ref{thm4} yields
$u\equiv v$ near $(0, \bar y)$. 
 Theorem \ref{thm8} follows in this case 
in view of Theorem \ref{thm6}.   

In Case 2, we still have (\ref{kkk}) and (\ref{hhh}) for small
$\epsilon>0$.  As explained in the proof of Theorem \ref{thm7},
we still have $\tau_s<1$ in $\Omega^+$.
Thus we still have
$L\tau\le 0<L(s+s^{\frac 32})$ in $\Omega^+_\epsilon$ and,
for some $c>0$, 
$\tau\ge c(s+s^{\frac 32})$ on $\partial \Omega^+_\epsilon$,
where
$$
L=\partial_{ss}+\Delta _y 
+(2-\tau_s)\frac{u_{tt}}{u_t}\partial_s-\eta\cdot \nabla_y\tau.
$$
Theorem \ref{thm8} in this case follows
as in the proof of Theorem \ref{thm7}.

\vskip 5pt
\hfill $\Box$
\vskip 5pt

Finally we include  the following 
result.

\begin{thm} Let  $u$ be a $C^\infty$ function
in the unit ball $B_1$ in $\Bbb R^n$, $n\ge 1$, satisfying
\begin{equation}
\Delta u(x)=V(x)u(x),\qquad x\in B_1,
\label{c1}
\end{equation}
where $V\in C^1(B_1\setminus\{0\})$ satisfies,
 in polar coordinates $(r,\theta)$, $\theta\in \Bbb S^{n-1}$,
\begin{equation}
(r^2V)_r:=\frac {\partial }{\partial r}(r^2V)\ge 0,
\qquad \mbox{in}\ B_1\setminus\{0\}.
\label{c20}
\end{equation}
Assume that $u$ vanishes of infinite order at the origin,
i.e.
\begin{equation}
\partial^{\alpha}u(0)=0\ \ \mbox{for all multi-index}\
\alpha=(\alpha_1, \cdots, \alpha_n),\ \alpha_i\ge 0.
\label{c30}
\end{equation}
Then $u\equiv 0$ in $B_1$.
\label{thmc1}
\end{thm}

\noindent{\bf Proof.}\  We make use of ideas in
Agmon and Nirenberg \cite{AN}.  
Using polar coordinates $(r,\theta)$, 
equation (\ref{c1}) takes the form
\begin{equation}
r^2 u_{rr}+(n-1)r u_r+\Delta_\theta u=r^2Vu.
\label{c2}
\end{equation}
Set
$$
r=e^s.
$$
Then
$$
u_s=u_r e^s, \ \ u_{ss}=u_{rr}e^{2s}+u_r e^s=
r^2 u_{rr}+ru_r,
$$
and 
(\ref{c2}) takes the form
$$
u_{ss}+(n-2) u_s +\Delta_\theta u=
r^2Vu,\qquad 
(s,\theta)\in (-\infty, 0)\times  \Bbb S^{n-1}.
$$
Because of (\ref{c30}),
\begin{equation}
\lim_{s\to -\infty}\max_{\theta\in \Bbb S^{n-1}}
\sum_{i=0}^2
(|\partial^i_su(s,\theta)|+|\partial^i_\theta
u(s.\theta)|)e^{b s}
=0,\qquad \forall\ b<0.
\label{c3}
\end{equation}
Set 
$$
u=e^{as}v\ \mbox{with}\
a=-\frac{n-2}2.
$$
Since
$$
u_s=e^{as} (v_s+av), \quad
u_{ss}=e^{as}(v_{ss}+2a v_s +a^2v),
$$
$v$ satisfies
\begin{equation}
v_{ss}+\Delta_\theta v=mv,\qquad
\mbox{in}\ (-\infty, 0)\times \Bbb S^{n-1},
\label{c4}
\end{equation}
where
$$
m:= (\frac {n-2}2)^2 +r^2V.
$$

Consider
$$
\rho(s):= \int_{\Bbb S^{n-1} }v^2(s, \theta)d\theta.
$$
We will prove
\begin{lem}
\begin{equation}
\frac{d^2}{ds^2} \log \rho(s)\ge 0\ \
\mbox{whenever}\ \rho(s)>0.
\label{c6}
\end{equation}
\label{lemc2}
\end{lem}

\noindent{\bf Proof.}\
By computation,
$$
\rho_s=2\int_{ \Bbb S^{n-1} }vv_sd\theta,\quad
\rho_{ss}=2\int_{ \Bbb S^{n-1} }
(v_s^2+ vv_{ss})d\theta.
$$
So
\begin{equation}
\rho_{ss}=2\int_{ \Bbb S^{n-1} }
v_s^2 +2 \int_{ \Bbb S^{n-1} }
v(-\Delta_\theta v+mv)d\theta
=2\int_{ \Bbb S^{n-1} }[ v_s^2+|\nabla_\theta v|^2 +mv^2]d\theta.
\label{c71}
\end{equation}
Next, by the Schwartz inequality,
\begin{equation}
\frac{ \rho_s^2}\rho =
\frac{ 4  (\int_{ \Bbb S^{n-1} }
vv_s d\theta)^2  }
{ \int_{ \Bbb S^{n-1} }
v^2 d\theta}\le
4 \int_{ \Bbb S^{n-1} }
v_s^2 d\theta.
\label{c72}
\end{equation}

Multiplying (\ref{c4}) by $2v_s$ and integrating in $s$ from $-\infty$
to $0$, and integrating over $\Bbb S^{n-1}$, we find, using
Green$'$s theorem,
$$
\int_{ \Bbb S^{n-1} }
v_s^2 -2
\int_{-\infty}^s \int_{ \Bbb S^{n-1} }
\nabla_\theta v_s \cdot \nabla _\theta v 
=\int_{-\infty}^s \int_{ \Bbb S^{n-1} }
2mvv_s,
$$
i.e.
%\begin{eqnarray*}
$$
\int_{ \Bbb S^{n-1} }
 v_s^2= \int_{ \Bbb S^{n-1} }
 |\nabla_\theta v|^2 +
 \int_{-\infty}^s \int_{ \Bbb S^{n-1} }
 2mvv_s
 =  \int_{ \Bbb S^{n-1} }
  |\nabla_\theta v|^2+  \int_{ \Bbb S^{n-1} }
  mv^2-  \int_{-\infty}^s \int_{ \Bbb S^{n-1} }
  m_sv^2.
 $$
 %\end{eqnarray*}
We know from  (\ref{c20}) 
that $m_s\ge 0$.  Thus
\begin{equation}
\int_{ \Bbb S^{n-1} }
 v_s^2\le
 \int_{ \Bbb S^{n-1} }
 [ |\nabla_\theta v|^2 +
  mv^2 ], \qquad s\in (-\infty, 0).
  \label{c8}
  \end{equation}
We deduce from (\ref{c71}), (\ref{c72}) and (\ref{c8}) that
$$
\rho_{ss}\ge \frac{ \rho_s^2 }\rho,
\qquad \mbox{whenever}\ \rho(s)>0,
$$
which is equivalent to (\ref{c6}).  Lemma \ref{lemc2}
is established.

\vskip 5pt
\hfill $\Box$
\vskip 5pt

To prove Theorem \ref{thmc1}, we only need to show that $\rho\equiv 0$.
Suppose $\rho(s)>0$
for some $\bar s\in (-\infty, 0)$.
By  (\ref{c6}), $\log \rho$ is convex in any open interval where
$\rho>0$.  So, for any interval $(-T, \bar s)$ where $\rho $ is positive,
we have
$$
\log \rho(s)\ge \log \rho(\bar s)
+\frac {d}{ds}  \log \rho(\bar s) (s-\bar s),
\qquad\forall\ -T<s<\bar s.
$$
It follows from the above that $\rho(s)>0$ for all $-\infty<s<\bar s$ and,
for some constant $C_1, C_2>0$,
$$
\rho(s)\ge C_1 e^{ C_2 s}, \qquad \forall\ -\infty<s< \bar s,
$$
which violates (\ref{c3}).  Theorem \ref{thmc1} is established.

\vskip 5pt
\hfill $\Box$
\vskip 5pt

\section{Appendix}

Let ${\cal S}^{n\times n}$ denote the set of
of real  $n\times n$ symmetric matrices,
and let $O(n)$ denote the set of real  $n\times n$
orthogonal matrices.   
For $N\in  {\cal S}^{n\times n}$, we use
$\displaystyle{
|N|:=\sqrt{\sum_{0\le k,l\le n-1}|N_{kl}|^2}
}
$
to denote the norm of $N$.
\begin{lem}  
Let $G$ be a $C^3$ function defined on 
${\cal S}^{n\times n}$ satisfying
$$
G(O^{-1}NO)=G(N),\qquad
\forall\ N\in {\cal S}^{n\times n}, \ \forall\
O\in O(n).
$$
Then, for some constant $C$ depending only on $n$ and $G$,
$$
\bigg|\sum_{\alpha=1}^{n-1}
\frac{\partial G}{ \partial N_{0\alpha} }(N)
N_{0\alpha}\bigg|
\le C\sum_{ \beta=1 }^{n-1}|N_{0\beta}|^2,\
\forall\ N\in {\cal S}^{n\times n},
\ |N|\le 1.
$$
\label{lema-1}
\end{lem}

\noindent{\bf Proof.} 
Let $\overline N$ denote elements in ${\cal S}^{n\times n}$
satisfying
$$
\overline N_{0\alpha}=\overline N_{\alpha 0}=0,\qquad
1\le \alpha\le n-1,
$$
and let $e$  denote elements in ${\cal S}^{n\times n}$
satisfying
$$
e_{00}=e_{\alpha\beta}=0,\qquad 1\le \alpha,\beta\le n-1.
$$
Consider the following function of $e$:
$$
h(e)=\frac d{dt} G(\overline N+te)\bigg|_{t=1}.
$$
For $N=\overline N+e$,
$$
h(e)= 2\sum_{\alpha=1}^{n-1}
\frac{\partial G}{ \partial N_{0\alpha} }(N)
N_{0\alpha}.
$$
Clearly
$$
h(0)=0.
$$
For $O=diag(-1, 1, \cdots, 1)$,
$$
O(\overline N+te)O=\overline N-te.
$$
So
$$
h(e)\equiv h(-e).
$$
Consequently
$$
\nabla h(0)=0.
$$
Since $h$ is a $C^2$ function, we obtain
$$
|h(e)|\le C|e|^2.
$$
Lemma \ref{lema-1} is established.

\vskip 5pt
\hfill $\Box$

\end{document}